\documentclass[11pt]{article}\hfuzz=10pt  \topmargin=-0.5cm\hfuzz=10pt  \oddsidemargin=0.1cm\textheight=220mm \textwidth=16.5cm
\usepackage{amssymb}
\usepackage{epsfig}
\newcommand{\comm}[1]{ }

\def\xxxonly{\comm}
\def\xxxonly{{ }}
\def\noxxx{\comm}

\newtheorem{theorem}{Theorem}
\newtheorem{lemma}{Lemma}
\newtheorem{condition}{Condition}
\newtheorem{remark}{Remark}
  
 \def\e{\varepsilon}
\def\tt{\theta}
\def\defi{\stackrel{{\scriptscriptstyle \Delta}}{=}}

\def\a{\alpha}
\def\d{\delta}

\def\w{\widehat}
\def\Ind{{\,\rm Ind\,}}

\def\essinf{\mathop{\rm ess\, inf}}

\def\R{{\bf R}}

\def\L{L}

\def\g{\gamma}

\def\W{{\cal W}^*}
\def\ww{\widetilde}

\def\t{\theta}
\def\oo{\bar}

\def\p{\partial}

\def\V{{\cal V}}

\def\M{{\cal M}}

\def\L{{\cal L}}

\newcommand{\be}{\begin{equation}}
\newcommand{\ee}{\end{equation}}
\newcommand{\bd}{\begin{displaymath}}
\newcommand{\ed}{\end{displaymath}}
\newcommand{\ba}{\begin{array}{ll}}
\newcommand{\ea}{\end{array}}
\newcommand{\baa}{\begin{eqnarray}}
\newcommand{\eaa}{\end{eqnarray}}
\newcommand{\baaa}{\begin{eqnarray*}}
\newcommand{\eaaa}{\end{eqnarray*}}
\font\sm=cmr10

\def\Ind{{\mathbb{I}}}
\def\L{{\cal L}}

\def\t{\tau}

\def\tt{\theta}

\def\W{{\cal W}}

\def\R{{\bf R}}

\def\t{\theta}
\def\M{M}
\date{{\xxxonly{Submitted September 7, 2016. Revised January 14, 2019 }}} 
\title
{On recovering parabolic diffusions from their time-averages}
\author{
Nikolai Dokuchaev}
\begin{document}
\maketitle
\let\thefootnote\relax
\footnote{
This is a an extended version of  an  article accepted for publication in {\em Calculus of Variations and Partial Differential Equations} \index{Accepted 9.12.2018}  and available online at https://doi.org/10.1007/s00526-018-1464-1.
The second example in Section 4 here  has been excluded form the journal version. 
 } 
\noxxx{\let\thefootnote\relax\footnote{The author is with  School of Electrical Engineering, Computing and Mathematical Sciences, Curtin
University,   GPO Box U1987, Perth, 6845 Western Australia. Email
N.Dokuchaev@curtin.edu.au. Ph. 61 8 92663144}
}
\begin{abstract}
The paper study a possibility to recover a parabolic diffusion
from its time-average when the values at the initial time are  unknown.  This problem
 can be reformulated as a new boundary value problem where a Cauchy condition is replaced by
 a prescribed time-average of the solution. It is shown that this new problem is well-posed in certain classes of solutions.
  The paper establishes existence, uniqueness,
 and  a regularity of the solution for this new problem and its modifications, including problems with singled out
 terminal values.
\end{abstract}
\section{Introduction}
Parabolic diffusion equations have fundamental significance for
natural and social sciences, and various boundary value problems for
them were widely studied including   inverse and ill-posed problems;
see examples in Miller (1973),  Tikhonov and Arsenin (1977), Glasko
(1984), Prilepko {\it et al} (1984),
  Beck (1985),   Showalter (1985),  Clark and  Oppenheimer  (1994),   Seidman (1996),
H\'ao (1998), Li {\em et al} (2009), Triet  {\em et al} (2013), Tuan and Trong (2011), Tuan and Trong (2014), Hao  (1998), Bourgeois and Dardé  (2010),
H\'ao and Oanh (2017), and the references therein.

According to Hadamard criterion, a boundary value problem is
well-posed if it features existence and uniqueness of the solution as well as
 continuous dependence of the solution on the data. Otherwise, a problem is ill-posed.

\par
For  parabolic equations, it is commonly recognized that the choice of the time
where the Cauchy condition is imposed  defines if a problem is well-posed
or ill-posed.
  A classical example is the
heat equation \baaa u'_t(x,t)=u_{xx}^{''}(x,t), \quad
t\in[0,T]. \label{sample}\eaaa The problem for this equation with
the Cauchy condition $u(x,0)\equiv \mu(x)$ at the initial time $t=0$ is well-posed in usual
classes of solutions.  In contrast, the problem with the Cauchy
condition  $u(x,T)\equiv \mu(x)$ at the terminal time $t=T$  is ill-posed. This means  that a
prescribed profile of temperature at time $t=T$ cannot be achieved
via an appropriate selection of the initial temperature. Respectively, the
initial temperature profile cannot be recovered from the observed temperature at the terminal time.
In particular, the process $u$ is not robust with respect to small deviations of its terminal profile $u(\cdot,T)$.
This makes this problem ill-posed, despite
the fact that solvability and uniqueness still can be achieved for
some very smooth analytical boundary data or for special selection
of the domains; see e.g. Miranker (1961), Dokuchaev (2007).
\par
It appears that there are boundary value problems that do not fit the
dichotomy of the classical forward/backward well-posedness. For instance,  the problems for forward heat equations are well-posed  with non-local in time conditions  connecting the values  at different times
such as
\baaa
u(x,0)-k u(x,T)=\mu(x) \quad \hbox{or}\quad u(x,0)+\int_0^Tw(t)u(x,t)dt=\mu(x),
\eaaa for  given $k\in \R$ and given functions $\mu$, $w$.  Some results for parabolic equations and stochastic PDEs
with these non-local conditions replacing the Cauchy condition  were obtained  in Dokuchaev (2004,2008,2011,2015). In these conditions, the singled out $u(\cdot,0)$
helped to counterbalance  the presence of the future values, given 
some restrictions on $k$ and $w$.

\par The present paper further extends the setting with
 mixed in time conditions.
The paper investigates solutions
$u(x,t)$ of the forward parabolic equations
with some new conditions, such as
\baaa
\int_0^T  u(x,t)dt=\mu(x)\quad \hbox{or}\quad k_1 u(x,T)+k_2\int_0^T  u(x,t)dt=\mu(x),
\eaaa
replacing a well-posed Cauchy condition $u(x,0)=\mu(x)$,  for 
 a given terminal time $T>0$, 
 a given function $\mu$, and given $k_i\in\R$.
A crucial difference with the setting from Dokuchaev (2015) is that the setting of the present paper does not  require that the initial value $u(\cdot,0)$ is singled out;  instead,  the initial  value $u(\cdot,0)$  is  presented as $u(\cdot,t)dt$  at $t=0$ only, i.e. under the integral, with a infinitively small weight at $t=0$.
 Moreover, the present paper allows a setting with  $k_1\neq 0$, i.e.  where only the terminal value $u(\cdot,T)$ is singled out. This is different from the  quasi-boundary value (QBV) method  used for recovery of initial conditions for the heat equations, where the boundary condition  $u(x,T)+\e u(x,0)=\mu(x)$ with small $\e>0$ is considered as a replacement for the ill-posed final condition $u(x,T)=\mu(x)$; see, e.g. Showalter (1985),  Clark and  Oppenheimer  (1994), Seidman (1996), Triet  {\em et al} (2013), Triet and Phong (2016).   A related but different setting  with observable  spatial integrals of the solutions for
parabolic equations was considered in  H\'ao and Oanh (2017). Li {\em et al} (2009) considered a related but different again setting with  solutions of  parabolic equations observable on certain subdomains.

 Formally, the new problems introduced in the present  with time averaging
do not fit the framework given by the classical theory of
well-posedness for parabolic equations based on the correct
selection of the time for a Cauchy condition. However, we found  that these
new problems are well-posed for  $\mu\in H^2$, i.e. if the second partial derivatives of $\mu$ are square integrable  (Theorem \ref{ThM}). This can be interpreted as an existence of a
diffusion with a prescribed average over a time interval. In addition, this can be interpreted as solvability of
the following inverse problem:  given $\int_0^T u(x,t)dt$ for all $x\in D$,
recover the entire process $u(x,t)|_{D\times [0,T]}$. It is shown below that
this problem is well-posed. This is an interesting result, because it is known that, for any $c>0$,  the knowledge of  values $u|_{D\times [c,T]}$
does not ensure restoring of the values $u|_{D\times [0,c)}$; this problem would be ill-posed.

This result can be applied, for example, to reduce the costs of data processing
 for the analysis of the dynamics of heat propagation: it suffices to collect, store, and transmit, only time averages of
 temperatures rather then the entire history.

The rest of the work is organized as follows. In Section \ref{SecS}, we introduce
boundary value problem with averaging over time.
In Section \ref{SecM}, we present the main result and its proof (Theorem \ref{ThM}), and we  discuss the properties of solutions of the suggested boundary value problems. A numerical example is given in Section \ref{SecN}.
 \section{Problem setting}\label{SecS}
Let  $D\subset \R^n$ be an open bounded connected domain with $C^2$ - smooth
boundary $\p D$, and let $T>0$ be a fixed number. We consider the boundary value problems
 \baa
&&\frac{\p u}{\p t}=A u+\varphi \quad \hbox{for}\quad (x,t)\in D\times (0,T),\label{eq}\\
&&\hphantom{x}u(x,t)=0 \quad \hbox{for}\quad (x,t)\in\p D\times
(0,T) \label{dD},\\
&&\kappa u(x,T)+\int_0^Tw(t)u(x,t)dt=\mu(x)\quad \hbox{for}\quad x\in D.\label{ppp}\eaa
Here  $\kappa\in \R$ and a function $w(t)$ are  given,
$$
A u\defi \sum_{i=1}^n  \frac{\p }{\p x_i}\left(\sum_{j=1}^n a_{ij}(x)\frac{\p u}{\p x_j}(x)\right)+a_0(x,t)u(x).
$$
 The functions $a_{ij}(x): D\to \R$ and $a_0(x): D\to
 \R$ are continuous and bounded, and there exist
 continuous bounded derivatives $\p a_{ij}(x,t)/\p x_i$, $i,j=1,...,n$.  In
addition, we assume that the matrix $a=\{a_{ij}\}$ is symmetric  and
$y^\top a(x)y\ge \d |y|^2$ for all $x\in D$ and $y\in\R^n$, where $\d>0$ is
a constant. The function $\varphi(x,t):  D\times(0,T)  \to \R$
is measurable and square integrable. Conditions (\ref{eq})-(\ref{dD}) describe a diffusion process in domain $D$.
\par
We consider problem (\ref{eq})-(\ref{ppp}) assuming that the coefficients of $A$ and the inputs  $\mu$ and $\varphi$ are known, and that the initial value
$u(\cdot,0)$  is unknown.

\par
 If $\kappa\neq 0$ and $w\equiv 0$, then problem (\ref{eq})-(\ref{ppp}) is ill-posed, with a  Cauchy condition $u(x,T)=\mu(x)$.
To exclude this case, we assume, up to the end of this paper, that the following condition holds.
\begin{condition} In (\ref{ppp}),  $\kappa\ge 0$, and the function $w$ is bounded and such that
\baaa
\quad  w(t)\ge 0\quad\hbox{for a.e.}\quad  t\in[0,T].
\eaaa
In addition, there exists $T_1\in(0,T]$ such that
$\essinf_{t\in[0,T_1]}w(t)>0$. 
\end{condition}

\subsubsection*{Some special cases}
\begin{enumerate}
\item
If $\kappa=0$ and $w(t)\equiv 1$, then condition (\ref{ppp}) becomes
\baa
\int_0^Tu(x,t)dt=\mu(x)\quad \hbox{for}\quad x\in D.\label{aver}\eaa
Problem (\ref{eq})-(\ref{dD}),(\ref{aver}) can be considered as a problem of recovering $u$
from its time-average $\int_0^Tu(x,t)dt$.
\item
If  $\kappa=1$, and $w(t)\equiv \Ind_{[0,\e]}(t)$ , then condition (\ref{ppp}) becomes
\baa
u(x,T)+\int_0^\e u(x,t)dt=\mu(x)\quad \hbox{for}\quad x\in D.\label{ill}
\eaa
With a small $\e>0$, solution of problem  (\ref{eq})-(\ref{dD}),(\ref{ill})  can  be considered as a variation of the quasi-boundary-value method
 for solution of backward equation, where an ill-posed condition $u(x,T)=\mu(x)$ is replaced by condition (\ref{ill}); see, e.g. Showalter (1985),  Clark and  Oppenheimer  (1994).\index{Clark,Sho} Seidman (1996),   Triet  {\em et al} (2013).

\end{enumerate}
Here $\Ind$ denotes the indicator function.

Some mild restrictions will be imposed on the choice of $\varphi$ for the case where $\kappa\neq 0$:
it will be required that $\varphi(\cdot,t)$ features some reqularity in $t\in[\t,T]$ for some $\t\in [0,T)$
that can be arbitrarily close to $T$.
\subsubsection*{Spaces and classes of functions}
For a Banach space $X$,
we denote  the norm by $\|\cdot\|_{ X}$. For a Hilbert space $X$,
we denote  the inner product by $(\cdot,\cdot)_{ X}$.

We denote
by $W_2^m(D)$   the
standard Sobolev spaces of functions that belong to $L_2(D)$ together with their generalized derivatives of $m$th order.
We denote by $\stackrel{0}{W_2^1}(D)$ the closure in the
${W}_2^1(D)$-norm of the set of all continuously
differentiable functions $u:D\to\R$ such that  $u|_{\p D}\equiv 0$; this is also a Hilbert space.
  \par
Let $H^0\defi L_2(D)$ and  $H^1\defi \stackrel{0}{W_2^1}(D)$.
\par
 Let $H^{-1}$ be the dual space to $H^{1}$, with the
norm $\| \,\cdot\,\| _{H^{-1}}$ such that if $u \in H^{0}$ then $\|
u\|_{ H^{-1}}$ is the supremum of $(u,v)_{H^0}$ over all $v \in H^1$
such that $\| v\|_{H^1} \le 1 $.
\par
Let $H^2$ be the subspace of $H^1$ consisting
of elements with a finite norm in $W_2^2(D)$; this is also a Hilbert space.
\par We
denote
 the Lebesgue measure and
 the $\sigma $-algebra of Lebesgue sets in $\R^n$
by $\oo\ell _{n}$ and $ {\oo{\cal B}}_{n}$, respectively.

\par
Introduce the spaces
$$
{\cal C}_k\defi C\left([0,T]; H^k\right),\quad \W^{k}\defi
L^{2}\bigl([ 0,T ],\oo{\cal B}_1, \oo\ell_{1};  H^{k}\bigr), \quad  k=-1,0,1,2,
$$
and the spaces
$$
\V^k\defi \W^{k}\cap {\cal C}_{k-1},\quad k=1,2,
$$
with the  norm $ \| u\|_{\V^k} \defi \| u\| _{{\W}^k} +\|
u\| _{{\cal C}_{k-1}}. $

For $\t\in [0,T)$, we introduce a space $\W_\t^0$ of functions $\varphi\in\W^0$ such that $\varphi(\cdot,t)=\oo\varphi+\int_\t^t\w\varphi(\cdot,s)ds$ for $t\in[\t,T]$ for
some $\oo\varphi\in H^0$ and $\w\varphi\in L_1([\t,T];H^0)$, with the norm
 \baaa\| \varphi\| _{\W^0_\t} \defi \| \varphi\| _{{\W}^0} +\|\oo\varphi\|_{H^0}+ \int_\t^T\|\w\varphi(\cdot,t)\|_{H^0}dt. \eaaa
In particular, $\varphi(\cdot,t)$ is continuous in $H^0$ in $t\in(T-\t,T]$.
 We extend this definition  on the case where $\t=T$, assuming that $\W^0_T=\W^0=L_2(D\times[0,T])$.

 \par
As usual, we accept that equations (\ref{eq})-(\ref{dD}) are
satisfied for $u\in \V^1$ if,  for any $t\in[0,T]$,
\baa
\label{intur1} u(\cdot,t)=u(\cdot,0)+\int_0^t [A u(\cdot,s)+\varphi(\cdot,s)]ds.
\eaa  The equality here is assumed to be an equality in the space
$H^{-1}$. Condition (\ref{ppp}) is satisfied
as an equality in $H^0=L_2(D)$. The condition on $\p D$ is satisfied in the  sense that
$u(\cdot,t)\in H^1$ for a.e.   $t$.  Further, we have that $A u(\cdot,s)\in
H^{-1}$ for a.e. $s$ and the integral in (\ref{intur1}) is
defined as an element of $H^{-1}$. Hence equality (\ref{intur1}) holds in the sense of equality in $H^{-1}$.
\section{The result}\label{SecM}
\begin{theorem}
\label{ThM} Let $\t\in [0,T]$ be such that $\t=T$ if $\kappa=0$ and $\t<T$
 if $\kappa\neq 0$.
 For any  $\mu\in H^2$ and $\varphi\in\W^0_\t$, there exists a unique solution
$u\in\V^1$ of problem (\ref{eq})-(\ref{ppp}). Moreover, there exists $c>0$  such
that
\baa \|u\|_{\V^1}^2\le c\left(\|\mu\|_{H^2}^2+\|\varphi\|_{\W^0_\t}^2\right)\label{estp} \eaa for all
$\mu\in H^2$ and $\varphi\in\W^0_\t$. Here $c>0$ depends only on $n,T,D,\t$, $\kappa$, $w$, and
on the coefficients of  equation (\ref{eq}). \end{theorem}
\par
By Theorem \ref{ThM}, problem
(\ref{eq})-(\ref{ppp}) is well-posed in the sense of Hadamard for $\mu\in H^2$ and $\varphi\in\W^0_\t$.
\par
The proof of this theorem is given below; it is based on  construction of the solution $u$ for given $\mu$ and $\varphi$.
\subsection{Proofs}
Let us introduce  operators $\L :H^k\to \V^{k+1}$, $k=0,1$, and $L: \W^k\to \V^{k+2}$, $k=-1,0$,  such that $\L\xi+L\varphi=v$,
where $v$ is the solution in $\V$ of  problem (\ref{eq})-(\ref{dD}) with the Cauchy condition
\baa
u(\cdot,0)=\xi\label{uxi}.
\eaa These linear operators are continuous; see
e.g. Theorems III.4.1 and IV.9.1 in Ladyzhenskaja {\it et al} (1968) or Theorem III.3.2 in Ladyzhenskaya  (1985).

Let a linear operator $M_0: H^0\to H^1$ be defined
such that \baaa
(M_0 \xi)(x)=\int_0^T w(t) u(x,t)dt+\kappa u(x,T),\quad u=\L\xi\in\V^1.
\eaaa In other words, $u$ is the solution
of problem (\ref{eq})-(\ref{dD}) with the Cauchy condition $u(\cdot,0)=\xi\in H^0$ and with $\varphi=0$.

Further, let a linear operator $M: \W^0\to H^1$ be defined
such that \baaa
(M \varphi)(x)=\int_0^T w(t)u(x,t)dt+\kappa u(x,T),\quad u=L\varphi\in\V_1.\eaaa In other words, $u$ is the solution
of problem (\ref{eq})-(\ref{dD}) with this $\varphi$ and with the Cauchy condition $u(\cdot,0)=0$.

In these notations, $\mu=M_0 u(\cdot,0)+M\varphi$ for  a solution $u$ of problem  (\ref{eq})-(\ref{dD}).
\begin{lemma}
\label{lemma0}  The linear operator $M_0: H^0\to H^2$  is a continuous bijection; in particular, the inverse operator $M_0^{-1}:H^2\to H^0$ is also continuous. Their norms depends only on $n,T,D,\t$, $\kappa$, $w$, and
on the coefficients of  equation (\ref{eq}).
\end{lemma}

\par
\begin{remark}\label{remC}{\rm
It can be noted that the classical results for parabolic equations imply that
the operators $M_0: H^k\to H^{k+1}$,  $k=0,1$, and  $M: \W^0\to H^{2}$,  are  continuous for $\kappa=0$, and
the operators $M_0: H^k\to H^{k}$,  $k=0,1$, and  $M: \W^0\to H^{1}$,  are  continuous for $\kappa > 0$; see Theorems III.4.1 and IV.9.1 in Ladyzhenskaja {\it et al} (1968)
or Theorem III.3.2 in Ladyzhenskaya  (1985). The continuity of the operator   $M_0: H^0\to H^2$ claimed in Lemma \ref{lemma0} requires a proof that is given below.
}\end{remark}

{\em Proof of Lemma \ref{lemma0}}.
It is known that there exists
an orthogonal basis  $\{v_k\}_{k=1}^\infty$  in $H^0$, i.e. such that
 $$(v_k,v_m)_{H^0}=0,\quad  k\neq m,\quad \|v_k\|_{H^0}=1,$$
such that $v_k\in H^1$ for all $k$, and that
\baa
 Av_k=-\lambda_k v_k,\quad  v_k|_{\p D}=0,\label{EP}\eaa
for some $\lambda_k\in\R$, $\lambda_k\to +\infty$ as $k\to +\infty$; see e.g.  Ladyzhenskaya (1985), Chapter 3.4. In other words,
$\lambda_k$ and $v_k$ are the eigenvalues and the corresponding
 eigenfunctions  of the eigenvalue problem (\ref{EP}).

If $u\in\V^1$ is a solution of problem (\ref{eq})-(\ref{ppp}) with  $\varphi =0$, then $u(\cdot,0)\in H^0$
is uniquely defined; it follows from the definition of $\V^1$. Hence $\xi=u(\cdot,0)\in H^0$
is uniquely defined. Let $\xi$ and $\mu$ be expanded  as
\baaa
\xi=\sum_{k=1}^\infty \a_k v_k,\quad \mu=\sum_{k=1}^\infty \g_k v_k,
\eaaa
 where
$\{\a_k\}_{k=1}^\infty$ and $\{\g_k\}_{k=1}^\infty$ and square-summable real sequences.
 By the choice of $\xi$, we have that  $u=\L \xi$. Applying the Fourier method, we obtain that
\baa
u(x,t)=\sum_{k=1}^{\infty} \a_k e^{-\lambda_k t}v_k(x).
\label{sol}\eaa
\par
On the other hand,
\baaa
&&
\mu(x)=\sum_{k=1}^{\infty} \g_k v_k(x)= \int_0^Tw(t)u(x,t)dt+\kappa u(x,T)\\&&=
\sum_{k=1}^{\infty} \int_0^Tw(t)\a_k e^{-\lambda_k t} v_k(x)dt + \kappa\sum_{k=1}^{\infty} \a_k e^{-\lambda_k T} v_k(x)
\\&&=
\sum_{k=1}^{\infty} \zeta_k \a_k v_k(x),
\eaaa
where
\baaa &&\zeta_k =\int_0^Tw(t)e^{-\lambda_k t}dt+\kappa e^{-\lambda_kT}.
\eaaa
Therefore, the sequence $\{\a_k\}$ is uniquely defined as
\baa
\a_k =\g_k/\zeta_k,
\quad k=1,2,....
\label{ak}
\eaa
Remind that we had assumed that there exists $T_1>0$ such that $w_*\defi \inf_{t\in[0,T_1]} w(t)>0$ and that $\kappa\ge 0$.
In particular, this implies that $\zeta_k>0$ for all $k$. Moreover,
we have that
\baaa
\zeta_k \ge w_*\int_0^{T_1}e^{-\lambda_k t}dt+\kappa e^{-\lambda_kT}= w_*\frac{1-e^{-\lambda_kT_1}}{\lambda_k}
+\kappa e^{-\lambda_kT}.
\eaaa In addition, we have that
 \baaa
\zeta_k \le w_+\int_0^{T_1}e^{-\lambda_k t}dt+\kappa e^{-\lambda_kT}= w_+\frac{1-e^{-\lambda_kT_1}}{\lambda_k}
+\kappa e^{-\lambda_kT},
\eaaa where $w_+\defi \sup_{t\in[0,T_1]} w(t)$,\par
By the properties of  $A$, we have that
$\lambda_k\to +\infty$ as $k\to +\infty$, and that this sequence is non-decreasing. Hence there exists $m\ge 0$  such that $\lambda_m>0$; respectively,
 $\lambda_k>0$ for all $k\ge m$.

 Let
  \baaa
 &&c_1=\min\left[\zeta_1,...,\zeta_m,\,w_* \left(1-e^{-\lambda_m T_1}\right)\right],\qquad \\
 &&c_2=\max\left[\zeta_1,...,\zeta_m,\,w_+\left(1-e^{-\lambda_m T_1}\right) +\kappa\sup_{\lambda>0}\lambda e^{-\lambda T}\right].
 \eaaa
 Clearly, $0<c_1<c_2$ and
\baa
&c_1\le \lambda_k\zeta_k\le c_2,\quad &k\ge m,\nonumber\\
& c_1\le \zeta_k\le c_2,\quad &k<m.\label{zeta}\eaa
This can be rewritten as
 \baaa
&c_2^{-1}\lambda_k\le \zeta_k^{-1}\le c_1^{-1}\lambda_k,\quad &k\ge m,\\
& c_2^{-1}\le \zeta_k^{-1}\le c_1^{-1},\quad &k<m.\label{zetaest}\eaaa
It can be noted that estimate (\ref{zeta}) is crucial for the proof; this estimate defines regularisation
with $T_1$ is a parameter.

It follows that there exist some $C_1>0$ and $C_2>0$ such that
\baa
\sum_{k=1}^{\infty} \a_k^2\le C_1\sum_{k=1}^{\infty} \g_k^2 \lambda_k^2\le C_2 \sum_{k=1}^{\infty} \a_k^2.
\label{3}\eaa
We have that \baaa
A\mu=\sum_{k=1}^{\infty} \g_k A v_k(x) =-\sum_{k=1}^{\infty} \g_k \lambda_k v_k(x)
\eaaa
and
\baa
\|A\mu\|_{H^0}^2=\sum_{k=1}^{\infty} \g_k^2 \lambda_k^2,\quad \|\xi\|_{H^0}^2=\sum_{k=1}^{\infty} \a_k^2  <+\infty.
\label{Amuf}
\eaa
Hence (\ref{3}) can be rewritten as
\baa
\|\xi\|_{H^0}^2\le C_1\|A\mu\|_{H^0}^2\le C_2\|\xi\|_{H^0}^2.
\label{3n}
\eaa
\par
Suppose that $\mu\in H^2$. In this case,  $\|A \mu\|_{H^0}\le C \| \mu\|_{H^2}$,  for some $C>0$ that is
independent on $\mu$. Thus, (\ref{3n})  implies that the operator $M_0^{-1}:H^2\to H^0$ is continuous.

Let us prove  that the operator $M_0:H^0\to H^2$ is continuous.
 From the  classical estimates
 for parabolic equations, it
follows  that the operator $\L:H^0\to \V^1$ is continuous; see, e.g., Theorem IV.9.1 in Ladyzhenskaja {\em et al} 
(1968). By the definition of the operator  $M_0$, it follows that the operator $M_0:H^0\to H^0$ is continuous.

Further,  suppose that $\xi\in H^0$ and $\mu=M_0\xi$. Since
the operator $M_0:H^0\to H^0$ is continuous, we have that $\mu\in H^0$.
By (\ref{3n}),  $A\mu\in H^0$.
It follows that, for any $\lambda\in \R$,  we have that $h\defi A\mu+ \lambda \mu\in H^0$. By the properties of the elliptic equations,
it follows  that there exists $\lambda\in\R$ and $c=c(\lambda)>0$ such that \baa
\|\mu\|_{H^2}\le c\|h\|_{H^0}\le c(\|A\mu\|_{H^0}+\|\lambda \mu\|_{H^0});
\label{2ndfund}\eaa
see e.g. Theorem II.7.2 and Remark II.7.1 in Ladyzhenskaya (1975), or  Theorem III.9.2 and Theorem  III.10.1 in Ladyzhenskaya and Ural'ceva  (1968).
By (\ref{2ndfund}), we have that  \baa
\|\mu\|_{H^2}\le c_1(\|A\mu\|_{H^0}+\|\xi\|_{H^0})\le c_2\|\xi\|_{H^0}
\label{2ndfund2}\eaa
for some $c_i>0$ that are independent on $\xi$ and depend only on $n,T,D,\t$, $\kappa$, $w$, and
on the coefficients of  equation (\ref{eq}).
This completes the proof of Lemma \ref{lemma0}.

We now in the position to prove Theorem \ref{ThM}.

{\em Proof of Theorem \ref{ThM}.} Let us show first that the operator $\M:\W^0_\t\to H^2$ is continuous.
 As was mentioned in Remark \ref{remC}, the operator $M:\W^0\to H^2$ is continuous for $\kappa=0$; in this case, we can select $\t=T$ and $\W^0_\t=\W^0=L_2(D\times[0,T])$.

Let us show that the operator $\M:\W^0_\t\to H^2$ is continuous for the case where $\kappa\neq 0$. By the assumptions, $\t<T$  in this case  and $\varphi(\cdot,t)=\oo\varphi+\int_\t^t\w\varphi(\cdot,s)ds$ for $t\in[\t,T]$ for
some $\oo\varphi\in H^0$ and $\w\varphi\in L_1([\t,T];H^0)$.  Without a loss of generality, let us assume  that  $\kappa=1$ and $w(t)\equiv 0$, i.e. $\mu =\M \varphi=u(\cdot,T)$; it suffices  because  the boundary value problem is linear.

Let $v_k$ and $\lambda_k$ be such as defined in the  proof of Lemma \ref{lemma0}.

    Let $\mu$, $\varphi$, and $\w\varphi$, be expanded  as
\baaa \mu=\sum_{k=1}^\infty \g_k v_k,\quad
\varphi(\cdot,t)=\sum_{k=1}^\infty \phi_k(t) v_k,
\quad \oo\varphi=\sum_{k=1}^\infty \oo\phi_k v_k,\quad \w\varphi(\cdot,t)=\sum_{k=1}^\infty \w\phi_k(t) v_k.
\eaaa
Here $\{\g_k\}_{k=1}^\infty$ and $\{\oo\phi_k\}_{k=1}^\infty$  are square-summable real sequences,  the sequence $\{\phi_k(t)\}_{k=1}^\infty \subset L_2(0,T)$ and $\{\w\phi_k(t)\}_{k=1}^\infty \subset L_1(0,T)$
are such  that \baaa
\sum_{k=1}^\infty\int_0^T |\phi_k(t)|^2dt<+\infty,\quad \int_\t^T\left(\sum_{k=1}^\infty|\w\phi_k(t)|^2\right)^{1/2}dt<+\infty.
\eaaa

Applying the Fourier method for $u=L \varphi$, we obtain that
\baa
\mu(x)=\sum_{k=1}^{\infty} \g_k v_k(x)=  u(x,T)=
\sum_{k=1}^{\infty} v_k(x)\int_0^T\phi_k(t) e^{-\lambda_k (T-t)}dt\nonumber\\
=\sum_{k=1}^{\infty} v_k(x)(p_k+q_k),
\label{sol2}\eaa
where
\baaa
p_k=\int_0^\t\phi_k(t) e^{-\lambda_k (T-t)}dt,\qquad q_k=\int_\t^T\phi_k(t) e^{-\lambda_k (T-t)}dt
\eaaa
Clearly, \baaa
 |p_k|\le e^{-\lambda_k (T-\t)}\int_0^\t|\phi_k(t)| e^{-\lambda_k (\t-t)}dt\le
T^{1/2}e^{-\lambda_k (T-\t)} \|\phi_k\|_{L_2(0,T)}.
\eaaa
Further, we have that \baaa
\lambda_k q_k= -\int_\t^Te^{-\lambda_k (T-t)}\w\phi(t)dt+\phi_k(T) -\oo\phi_k e^{-\lambda_k (T-\t)}.
\eaaa
It follows that
\baaa
\sum_{k=1}^\infty \lambda_k^2p_k^2 +\sum_{k=1}^\infty \lambda_k^2q_k^2\le c\|\varphi\|^2_{\W^0_\t}
\eaaa
for some $c>0$ that does not depend on $\varphi$ and depends only on $n,T,D,\t$, $\kappa$, $w$, and
on the coefficients of  equation (\ref{eq}). Hence
\baaa
\|A\mu\|_{H^0}^2=\sum_{k=1}^\infty \lambda_k^2\g_k^2\le 2\sum_{k=1}^\infty \lambda_k^2p_k^2 +2\sum_{k=1}^\infty\lambda_k^2q_k^2
\le 2c\|\varphi\|^2_{\W^0_\t}.
\eaaa
 Similarly to (\ref{2ndfund})-(\ref{2ndfund2}), we obtain that  $\|\mu\|_{H^2}\le c\|A\mu\|_{H^0}$ for some $c>0$
that does not depend on $\varphi$ and depends only on $n,T,D,\t$, $\kappa$, $w$, and
on the coefficients of  equation (\ref{eq}). Hence the operator $\M:\W^0_\t \to H^2$ is continuous 
and its norm depends only on $n,T,D,\t$, $\kappa$, $w$, and
on the coefficients of  equation (\ref{eq}).

Further, it follows from the definitions of $\M_0$ and $\M$ that
\baaa
\mu=M_0\xi+M\varphi.
\eaaa
Since the operator $\M:\W^0_\t \to H^2$  and $\M_0^{-1}: H^2 \to H^0$ are  continuous, it follows that   $\M\varphi\in H^2$ and \baa
\xi=\M_0^{-1}(\mu-\M\varphi)\label{xi}
\eaa
is uniquely defined in $H^0$. Hence \baa
u=\L\xi+L\varphi=\L \M_0^{-1}(\mu-\M\varphi)+L\varphi.
\label{us}
\eaa
 is an unique  solution of problem (\ref{eq})-(\ref{ppp}) in $\V^1$.
  By the continuity of this and other operators in (\ref{us}), the
  desired estimate for $u$  follows.
This completes the proof of Theorem \ref{ThM}.
$\Box$

\begin{remark}{\rm
Equations (\ref{sol})--(\ref{ak})   provide a numerical method for calculating $\xi=M_0^{-1}\mu$. This and (\ref{us})
gives a  numerical method for
solution of problem
(\ref{eq})-(\ref{ppp}).}
\end{remark}
\subsection{On the properties of the solution}
The solutions of new problem (\ref{eq})-(\ref{ppp})   presented in Theorem \ref{ThM} have certain special  features described below.
\subsubsection*{Weaker regularity than for the classical problem}
It appears that  the solution of  new problem (\ref{eq})-(\ref{ppp})  has "weaker" smoothing properties than the solution of the classical
problem with standard initial Cauchy conditions. This can be seen from the fact that  problem
(\ref{eq})-(\ref{dD}),(\ref{uxi}) is solvable in $\V^2$ with a initial value $u(\cdot,0)\in H^1$ and with $\varphi\in \W^0$,
In addition, standard  problem
(\ref{eq})-(\ref{dD}),(\ref{uxi}) is solvable in $\V^1$ with  $u(\cdot,0)\in H^0$ and $\varphi\in \W^{-1}$.
   On the other hand,
new problem  (\ref{eq})-(\ref{ppp}) with $\mu\in H^2$  provides solution in $\V^1$ only, and does not allow $\varphi\in \W^{-1}\setminus \W^{0}$.
\subsubsection*{Non-preserving non-negativity}
For the classical problem (\ref{eq})-(\ref{dD}),(\ref{uxi})
with the standard Cauchy condition $u(x,0)=\xi(x)$, we have that if $\xi(x)\ge 0$ and $\varphi(x,t)\ge 0$ a.e. then $u(x,t)\ge 0$ a.e.
This is so-called Maximum Principle for parabolic equations;
see e.g. \cite{LSU}, Chapter III.7).

It  appears that  this does not hold for condition (\ref{ppp}):
a solution of problem (\ref{eq})-(\ref{ppp}) with non-negative functions $\mu$ and $\varphi$ is not necessarily non-negative.
It follows from the  Maximum Principle for parabolic equations that if $\xi(x)=u(x,0)\ge 0$ a.e. then $\mu(x)=(M_0\xi)(x)\ge $ a.e..
However, it may happen that the function $u(\cdot,0)=M_0^{-1}\mu$
can take negative values even if  $\mu(x)>0$ in all interior points of $D$.  This is because $\mu=M_0u(\cdot,0)$ actually represents a smoothing of $u(\cdot,0)$, and
this smoothing is capable of removing small negative deviations of $u(\cdot,0)$. This feature is illustrated by a numerical example in Section \ref{SecN} below.

\subsubsection*{A stability and robustness in respect to deviation of $\mu$ in $H^2$}
Let us discuss stability of the solution  implied by Theorem \ref{ThM}, or robustness
 in respect to deviation of $\mu$ in $H^2$.  Let us considered a family of functions
\baaa
\mu_{\d}(x)=\mu(x)+\d \eta(x),\quad \varphi_{\d}(x,t)=\varphi(x,t)+\d \psi(x,t), \quad  \d>0,
 \label{mudM} \eaaa where  $\eta\in H^2$ and $\psi\in\W^0_\t$ represent deviations.
 Let $u_\d$ be the corresponding solutions of problem (\ref{eq})-(\ref{ppp}). It follows from the linearity of the problem
 that
 \baaa
 \|u_0-u_\d\|_{\V^1}\le c \d\left(\|\eta\|_{H^2}^2+\|\psi\|_{\W^0_\t}^2\right), \eaaa
where $c>0$ is the same as in (\ref{estp}); this shows that the solution is robust with respect to deviations of inputs.

However, this robustness has its limitations  since
the norm $\|\eta\|_{H^2}$ can be large for non-smooth or frequently oscillating  $\eta$. For example, consider
$\eta(x)=\eta_\tt(x)=\sin(\tt x_1)\oo\eta(x)$, where $\tt>0$, $\oo\eta\in H^2$ is fixed and  $x_1$ is the first component of $x=(x_1,...,x_n)$.
In this case,  $|\eta_\tt(x)|\le |\oo\eta(x)|$ and $\|\eta_\tt\|_{H^2}\to +\infty $ as $\tt\to +\infty$  for a typical $\oo \eta$.
This feature is also illustrated by a numerical example in Section \ref{SecN} below.

\section{A numerical example}\label{SecN}
\xxxonly{\subsubsection*{An example  for $\mu$ defined by (\ref{aver})}}

Let us  consider  a numerical example for one-dimensional case where  $n=1$ and $D=(0,L)$. 
Let us consider a problem
\baa u'_t=u_{xx}''-q u,\quad  u|_{\p D}=0,\quad \int_0^Tu(x,t)dt=\mu(x),
\label{ex}\eaa
 where $q\ge 0$ is given.

To illustrate some robustness with respect to small deviations of $\mu$,  we considered a family of functions
\baa
\mu_{\d,\tt}(x)=\mu(x)+\d \eta_\tt(x),\quad \d>0,\quad \tt>0,
 \label{mud} \eaa where functions $\eta_\tt:D\to\R$ represent deviations and selected such that
the norm $\|\eta_\tt\|_{H^2}$ is increasing in  $\tt$    and that $\sup_{x}|\eta_\tt(x)|$ is bounded in $\tt$.

To solve the problem numerically, we calculated corresponding  truncated series
\baa
u_{\d,\tt,N}(x,0)=\sum_{k=1}^N\a_{k,\d,\tt} v_k(x).
\label{mud2}
\eaa using (\ref{sol}), (\ref{ak}) with $t=0$ and with corresponding $\a_k=\a_{k,\d,\tt}$.

For calcualtions, we have used $L=2\pi$, $q=0.0001$, $T=0.1$,
 $N=50$, and  $\tt=1,3$, and inputs
\baa
&&\mu(x)=x^{1/4}(L-x)|\sin(\pi x/L)|, \quad\nonumber\\&& \eta_\tt(x)=  x(L-x)\left(x-\frac{L}{3}\right)\left(x-\frac{2L}{3}\right)\sin (\tt x).
\label{mud1}\eaa
With this choice, the norms $\|d^2\eta_{\tt}(\cdot )/dx^2\|_{H^0}$ and   $\|\eta_\tt\|_{H^2}$ are increasing in  $\tt$.
  
  Some experiments with larger  $N=1000$ produced results that were almost indistinguishable from the results for $N=50$; we omit them here.

We have used MATLAB; the calculation for a standard PC  takes  less than  a second of CPU time, including calculation with larger $N>1000$.

\begin{figure}[ht]
\centerline{\psfig{figure=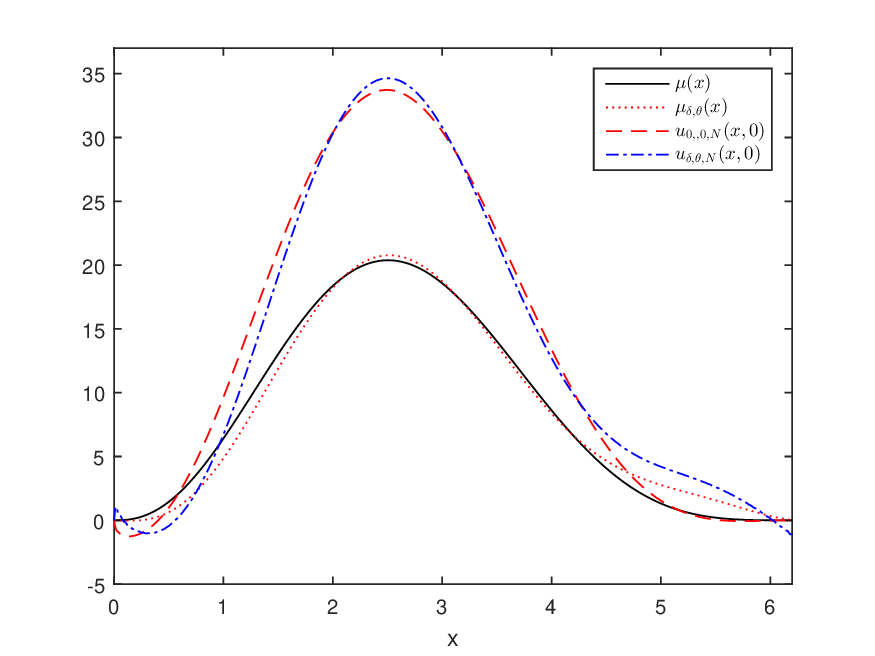,width=12cm,height=7cm}}
\centerline{\psfig{figure=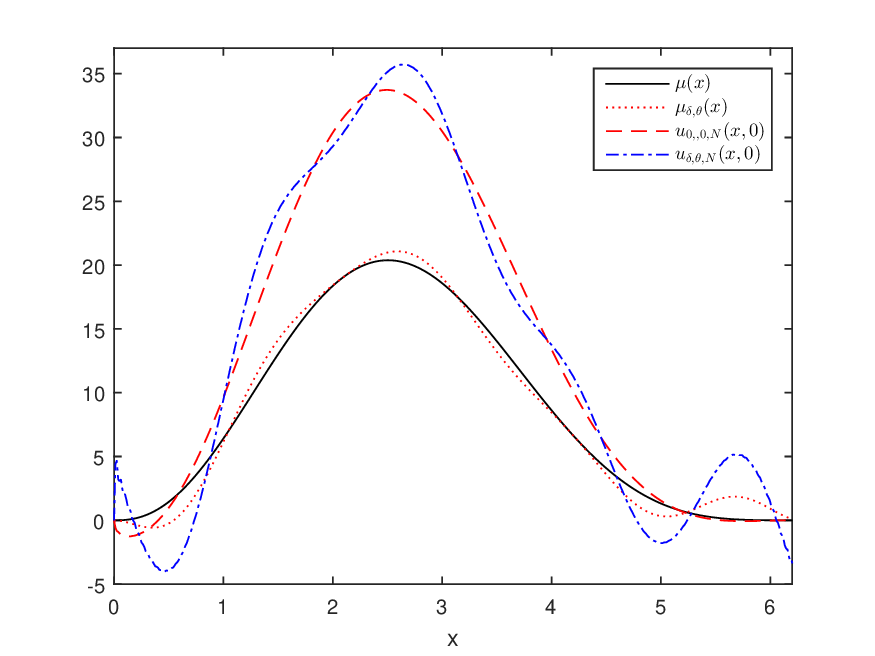,width=12cm,height=7cm}}
\caption[]{\sm The profiles of the time averages $\mu(x)$ and $\mu_{\d,\t}(x)$, and traces of the corresponding solutions $u_{0,0,N}(x,0)$
 and $u_{\d,\tt,N}(x,0)$ defined by (\ref{mud})-(\ref{mud2}) with $q=0.0001$, $T=0.1$, $\d=0.1$, $N=50$, $\tt=1$ (top) and $\tt=3$ (bottom).    }
\vspace{0cm}\label{f1}\end{figure}

 Figure \ref{f1} shows examples of  time averages $\mu$ and $\mu_{\d,\t}(\cdot)$,
and  corresponding profiles  $u_{\d,N,\t}(\cdot,0)$  recovered from the time averages via solution of problem (\ref{ex}) for $\d=0.1$
and for two choices $\t=1$ and $\t=3$.

Table \ref{simulation}  shows the relative  error
    \baaa
      E_{\d,N,\t}=\frac{\|u_{\d,N,\t}(\cdot,0)-u(\cdot,0)\|_{L_2(D)}}{\|u(\cdot,0)\|_{L_2(D)}}\label{err}\eaaa
of  recovery $u(x,0)$    calculated for a variety of $(\d,\t)$.

It can be seen from Figure \ref{f1} and Table \ref{simulation} that the solution is stable, i.e. it is robust  with respect
to small deviations of $mu$ in $H^2$. However, it can be also seen   that the magnitude of
 deviations of $u_{\d,\tt,N}(x,0)$ from $u_{0,0,N}(x,0)$ is larger for a larger  $\tt$.
 As was discussed in Section \ref{SecM}, this is consistent with   Theorem \ref{ThM}, because this theorem
 ensures robustness of the solutions with respect to deviations of $\mu$ that are small in $H^2$-norm. Respectively, deviations  that are small in  $H^0$-norm but large in $H^2$-norm may cause large deviations of solutions.

 Figure \ref{f1}  illustrates
 the comment in  Section \ref{SecM} pointing out on possibility to have non-negative solution of problem (\ref{eq})-(\ref{ppp}) for nonnegative $\mu$ and $\varphi$. The solution shown in Figure \ref{f1} have negative values, even  given that $\mu(x)>0$ for all $x\in D$.

\begin{table}[h]
\caption{Dependence of the relative error  $E_{\d,N,\t}$ on the input deviations.}
  \centering
\begin{tabular}{c|cccc}
  \hline
 &$\d=0.0001$&$\d=0.001$&$\d=0.05$&$\d=0.1$\\
  \hline
   $\t=0.05$  & 0.00002 & 0.00023 & 0.0113 & 0.0226\\
  $\t=0.1$  & 0.00004 & 0.00044 & 0.0218  & 0.0436\\
$\t=1$  & 0.00009 &  0.00087 & 0.0433  & 0.0866\\
 $\t=3$   &0.00014 & 0.0014& 0.0686 &  0.1372
\end{tabular}
\label{simulation}
\end{table}

\xxxonly{

\subsubsection*{An example for $\mu$ defined by (\ref{ill}) with applications to backward equations }
  By Theorem \ref{ThM},
$u(\cdot,0)$ can be restored from observation of  $\mu=\mu_\e$  for an arbitrarily small $\e>0$,
where $u$ is a solution of problem (\ref{eq})-(\ref{dD}),(\ref{ill}).
The following example illustrates a possibility to use this for the classical
problem of restoration of $u(\cdot,0)$ from $u(\cdot,T)$. For this problem,  $\mu=\mu_\e$ defined by (\ref{ill})
 is actually   unavailable for $\e>0$; instead, $u(\cdot,T)$ is available.
 Following the approach from Showalter (1985) and  Clark and  Oppenheimer  (1994),\index{Clark,Sho}
  we presume that the integral term in (\ref{ill}) is small, and we  accept $u(\cdot,T)$ as an approximation of $\mu_\e$.
 This leads to acceptance of
\baaa
u_\e(\cdot,0)\defi M_{\e,0}^{-1}u(\cdot,T)
\label{umu}\eaaa
 as an approximation of $u(\cdot,0)$, where $M_{\e,0}$ is defined as $M_0$ with
$\mu=\mu_\e$ defined by (\ref{ill}).
\par
We did some numerical experiments  to demonstrate potential applicability of this method.
Figure \ref{f2} demonstrates the results
for an example with $n=1$, $D=(0,L)$,  and with  the equation $u'_t=u_{xx}''-q u$, where $q>0$, $L>0$.
 In these  experiments, we first selected some profile $u(\cdot,0)$, then
 calculated  $u(\cdot,T)$ using the corresponding
  Green's function which is known for this toy forward equation; see e.g. \cite{But}, Chapter I.13. \index{p56rus}
  It can be noted that, for our experiment, it was sufficient to  use for the Green's function truncated  sin series with 50 terms. Further, for this
  $u(\cdot,T)$, we calculated $u_\e(\cdot,0)\defi\M_{\e,0}^{-1}u(\cdot,T)$
 using  equations (\ref{sol})--(\ref{ak}). Finally, we compared  $u_\e(\cdot,0)\defi\M_{\e,0}^{-1}u(\cdot,T)$ with true  $u(\cdot,0)$.

 More precisely, we used truncated series
\baa
u_{\e,N}(x,0)=\sum_{k=1}^N\a_{k,\e} v_k(x),\quad N>0,
\label{trunk}
\eaa
 as an approximation of the solution, where $\a_{k,\e}$ are defined by
(\ref{sol})--(\ref{ak}) applied for $w=w_\e$.

The limit case where $\e=0$ was not excluded; in this case,\baa
u_{0,N}(x,0)=\sum_{k=1}^N e^{\lambda_k T}g_{k} v_k(x)
\label{direct}\eaa is
a solution  based on straightforward truncation of the basis of eigenfunctions. Here $g_k\defi (u(\cdot,T),v_k)_{H^0}$.
For comparison purpose, we calculate this  solution as well.
\par
In addition,  we calculated an
estimate
 \baa
\ww u_{\e,N}(x,0)=\sum_{k=1}^N\frac{1}{\e+e^{-\lambda_k T}} g_{k} v_k(x).
\label{wwu}
\eaa
This estimate is implied by the quasi-boundary-value method  that suggests to replace  a ill-posed boundary condition $u(x,T)=f(x)$
 by a well-posed condition $\e u(x,0)+u(x,T)=f(x)$ such as in Showalter (1985),  Clark and  Oppenheimer  (1994).\index{Clark,Sho}

 Figure
\ref{f2} shows the results for  recovering
$u(x,0)=\Ind_{\{x>1.5\}}$  using our method with $\e=0.02$ and $N=18$. This figure shows $ u_{\e,N}(x,0)$ (our method), $\ww u_{\e,N}(x,0)$ (quasi-boundary-value method),
 and
$u_{0,N}(\cdot,0)$ (straightforward truncation (\ref{direct})).
Since $\e^{-1}\int_0^\e u(x,t)dt\approx u(x,0)$ in $L_2(D)$, it is natural to expect that the error for our  solution and  estimate (\ref{wwu}) implied by the quasi-boundary-value  method generate similar errors;  Figure
\ref{f2} shows that this holds for this example. In addition, it can be seen that these errors are  less than the error for the estimate defined (\ref{direct}).
It can be also noted that  $u_{0,N}(x,0)$ defined by (\ref{direct}) blows up for  $N\ge 19$.
Since analysis of the backward parabolic equations is not in  the focus of the present paper,
we leave the future research the questions of selection of $N$ and $\e$,  convergence analysis, and more precise comparison of different methods.

We used MATLAB and a standard PC; the calculation  takes  less than  a second of CPU time
for $N=1000$ in the setting of Figure \ref{f1}, and for $N=100$ in the setting of Figure \ref{f2}.


 \begin{figure}[ht]
\centerline{\psfig{figure=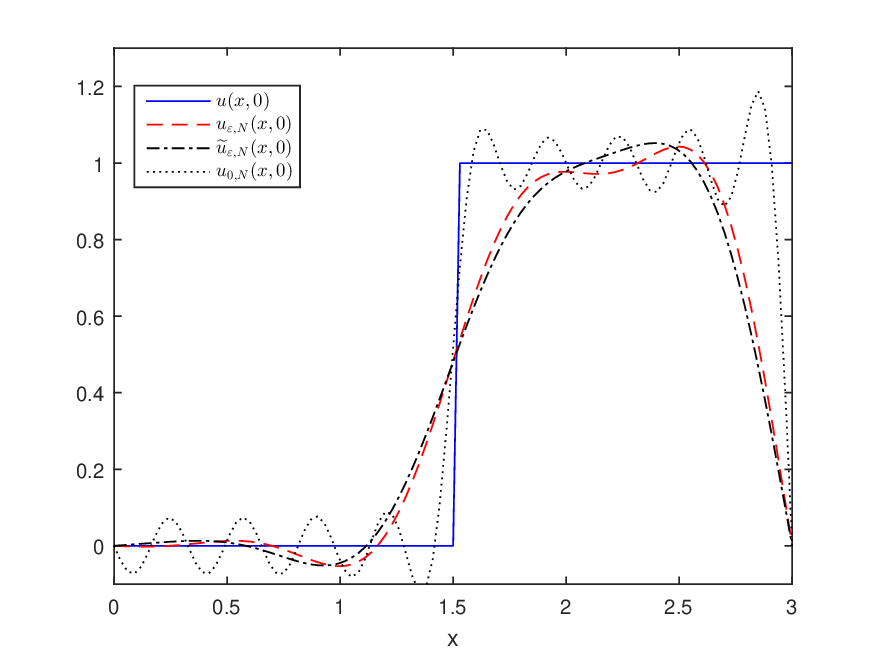,width=12cm,height=7cm}}
\caption[]{\sm
An initial profile $u(x,0)=\Ind_{\{x>1.5\}}$ and its estimates calculated for
$D=(0,3)$, $N=18$, $T=0.2$, and $\e=0.05$. Here $ u_{\e,N}(x,0)$
is estimate (\ref{umu}),  $\ww u_{\e,N}(x,0)$ is  estimate (\ref{wwu}),
$ u_{0,N}(x,0)$ is estimate (\ref{direct}).
 }
\vspace{0cm}\label{f2}\end{figure}
We used MATLAB and a standard PC; the calculation  takes  less than  a second of CPU time the calculation  takes  less than  a second of CPU time  for $N=100$ in the setting of Figure \ref{f2}. 
}

\section{Conclusion}\label{SecD}

The paper study a possibility to recover a parabolic diffusion
from its time-average for the case where the values at the initial time are  unknown.  This problem
 is reformulated as a new boundary value problem where a Cauchy condition is replaced by
 a condition involving the  time-average of the solution.   The paper establishes existence, uniqueness,
 and  a regularity of the solution for this new problem and its modifications, including problems with singled out
 terminal values (Theorem \ref{ThM}).
 This
Theorem \ref{ThM} can be applied, for example,  to the analysis of the evolution of temperature
in a domain $D$, with a fixed temperature on the boundary. The process $u(x,t)$ can be interpreted as
the temperature at a point $x\in D$ at time $t$. By Theorem
\ref{ThM}, it is possible  to recover the entire
evolution of the temperature in the domain if one knows the average temperature
over time interval $[0,T]$.

The suggested approach allows many modifications. An analog of Theorem \ref{ThM} can be obtained for the setting where problem (\ref{eq})--(\ref{ppp}) is considered for a known pair $(u(\cdot,0),\mu)$ and for
unknown $\varphi$ that has to be recovered. In this case, uniqueness of recovering $\varphi$ can be ensured via additional restrictions on its dependence on time; for example,
 it suffices to require that $\varphi(x, t)=\psi(t)v(x)$,
where $\psi $ is a known function,  and where  $v\in H^0$ is unknown and has to be recovered.

 It would be interesting to extend the result on the case where the operator $A$ is not necessarily symmetric and has coefficients depending on time.  We leave this for the future research.


\begin{thebibliography}{100}
\bibitem{B}
Beck, J.V. (1985). {\it
Inverse Heat Conduction.} John Wiley and Sons, Inc..

\noxxx{\bibitem{Bou}
Bourgeois, L and Dardé J. (2010)
A quasi-reversibility approach to solve the inverse obstacle problem, Inverse Problems and Imaging, vol.4, issue.3, pp.351-377, 2010.
}
\xxxonly{
\bibitem{But}
Butkovskiy, A. G. (1982).
{\em Green's Functions and Transfer Functions Handbook}, Halstead Press - John Wiley \& Sons, New York.
}
\bibitem{Clark}
 Clark G. W.,  Oppenheimer S. F. (1994)  Quasireversibility methods for non-well posed problems,
{\em Electronic Journal of Differential Equations },  no. 8, 1-9.

\bibitem{D04}
 Dokuchaev, N.G. (2004). Estimates for distances between first exit
times via parabolic equations in unbounded cylinders. {\it
Probability Theory and Related Fields}, {\bf 129} (2), 290 -
314.

\bibitem{D07} Dokuchaev, N. (2007). Parabolic equations with the second
order Cauchy conditions on the boundary. {\it  Journal of Physics A:
Mathematical and Theoretical}.
 {\bf 40}, pp. 12409--12413.
\bibitem{D08}
 Dokuchaev N. (2008). Parabolic Ito equations with mixed in time
conditions. {\it Stochastic Analysis and Applications} {\bf 26},
Iss. 3, 562--576.


\bibitem{D11}  Dokuchaev, N. (2011). On prescribed  change of
profile for solutions of parabolic equations.  {\it Journal of
Physics A: Mathematical and Theoretical} {\bf 44} 225204.

\bibitem{D15}
 Dokuchaev, N. (2015). On forward and backward SPDEs with non-local boundary
conditions. {\em Discrete and Continuous Dynamical Systems Series A (DCDS-A)}
{\bf 35}, No. 11, pp. 5335--5351

\bibitem{G} Glasko V. (1984). {\em Inverse problems of mathematical
physics}. American Institute of Physics. New York.
\bibitem{HaoB}
H\'ao, D.N.  (1998).  Methods for inverse heat conduction problems. Frankfurt/Main, Bern, New York, Paris: Peter Lang Verlag.
\bibitem{Hao}
Hao, D.N., Oanh, N.T.N. (2017). Determination of the initial condition in parabolic equations from integral observations. {\em  Inverse Problems
in Science and Engineering}  25 (8), 1138-1167.


\bibitem{L85}
Ladyzhenskaya, O. A. (1985). The boundary value problems of mathematical physics. Berlin etc., Springer-Verlag.
\bibitem{LSU}
Ladyzhenskaja, O.A.,  Solonnikov, V.A., and   Ural'ceva, N.N.
(1968). {\it Linear and quasi--linear equations of parabolic type.}
Providence, R.I.: American Mathematical Society.
\bibitem{LU68}
Ladyzhenskaja, O.A., and Ural'ceva, N.N., Linear and quasilinear elliptic equations, Academic Press, New York, 1968.
\bibitem{LiYam}
Li J., Yamamoto M., Zou J. (2009). Conditional stability and numerical reconstruction of initial
temperature. Commun. Pure Appl. Anal. 8, pp. 361?382.
\bibitem{M}
Miller, K. (1973). Stabilized quasireversibility and other nearly
best possible methods for non-well-posed problems. {\em In:
Symposium on Non-Well-Posed Problems and Logarithmic Convexity.
Lecture Notes in Math.} V. 316, Springer-Verlag, Berlin, pp.
161--176.
\bibitem{Mir} Miranker, W.L. (1961). A well posed problem for the
backward heat equation. {\em Proc. Amer. Math. Soc.} 12 (2), pp.
243-274.
\bibitem{P} Prilepko A.I., Orlovsky D.G.,
Vasin I.A. (1984). {\em Methods for Solving Inverse Problems in
Mathematical Physics}. Dekker, New York.
\bibitem{S}
Seidman, T.I. (1996). Optimal filtering for the backward heat
equation, {\em  SIAM J. Numer. Anal.} {\bf 33}, 162-170.
\bibitem{Sho}
Showalter, R.E.. (1985). Cauchy problem for hyper-parabolic partial differential
equations. In: Lakshmikantham, V. (ed.), Trends in the Theory and Practice of Non-Linear Analysis, Elsevier, North-Holland,  pp. 421-425.
\bibitem{T} Tikhonov, A.
N. and Arsenin, V. Y. (1977). {\it Solutions of Ill-posed Problems.}
W. H. Winston, Washington, D. C.
\bibitem{triet} Triet, L.M. and Phong, L.H. (2016). Regularization and error estimates for
asymmetric backward nonhomogenous heat equation in a ball. {\em Electronic
Journal of Differential Equations}, Vol. 2016, No. 256, 
pp. 1--12.
\bibitem{triet2}
Triet, M. L.; Quan, P. H.; Trong ,D. D.; Tuan, N. H.. (2013). A backward parabolic equation with
a time-dependent coefficient: Regularization and error estimates, J. Com. App. Math., No
237, pp. 432--441.
\bibitem{Tuan}
Tuan, N. H. and Trong, D. D. (2011). A simple regularization method for the ill-posed evolution equation. {\em Czechoslovak mathematical journal}  {\bf 61} (1), pp. 85-95.
\bibitem{Tuan2} Tuan N. H. and Trong, D. D . (2014). On a backward parabolic problem with local Lipschitz source J. Math. Anal. Appl. {\bf 414} 678?692.
\end{thebibliography}
\end{document}